\documentclass[12pt]{amsart}
\usepackage{amssymb} 
\usepackage{amsrefs}

\newtheorem{thm}{Theorem}[section]
\newtheorem{cor}[thm]{Corollary}

\theoremstyle{definition}

\newcommand{\Cl}{\mathrm{Cl}}
\newcommand{\rarr}{\rightarrow}
\newcommand{\Reg}{\mathrm{Reg}}

\newcommand{\RR}{{\mathbb R}}
\newcommand{\sub}{\subseteq}
\newcommand{\rdim}{\mathrm{dim}_\RR}
\newcommand{\CC}{{\mathbb C}}
\newcommand{\NN}{{\mathbb N}}

\renewcommand{\leq}{\leqslant}

\title{Subanalytic sets and complex analytic geometry}

\author[Peterzil]{Ya'acov~Peterzil}
\address{Department of Mathematics,
University of Haifa, Haifa, ISRAEL} \email{kobi@math.haifa.ac.il}

\author[Starchenko]{Sergei~Starchenko}
\address{Department of Mathematics, University of Notre Dame,
Notre Dame, IN 46556} \email{starchenko.1@nd.edu}

\date{September 21, 2004}
\dedicatory{Preliminary Report}%

\begin{document}
\maketitle

In a series of papers (\cite{ps1}, \cite{ps2}, \cite{ps3}) we
considered holomorphic manifolds and maps definable in  o-minimal
structures, over arbitrary real closed fields.  In this short note
we translate  some of these and more recent results into the
subanalytic setting. The theorems are of the following nature: We
consider a subanalytic subset $A$ of a complex analytic manifold
$M$ (when $M$ is viewed as a real manifold) and formulate
conditions under which $A$ is a complex analytic subset of $M$.\\

\noindent{\bf Remarks.} \\ 1. Even though we formulate the results
in the more familiar subanalytic language, they {\em all} remain
true when the subanalytic category is replaced by any
analytic-geometric category over $\RR$ (in the sense of van~den~Dries--Miller \cite{vddmil});
$\mathcal{C}_\mathrm{an,exp}$ is an example of such a category.\\
2. Because we prove most of the theorems  in a very general
o-minimal setting, most of the results yield certain uniformity in
parameters in the appropriate category.  \\

\noindent {\bf Notation and terminology.} For a real analytic
manifold $X$ and a subanalytic subset $A$ of $X$ we will denote by
$\rdim(A)$ the usual subanalytic dimension of $A$, i.e. the maximal $d\in \NN$
such that $A$ contains a $d$-dimensional $C^1$ submanifold of $X$.

We say that a set $S\subseteq \CC^n$ is {\em a subanalytic subset
of $\CC^n$ } if $S$ is a subanalytic  subset of $\RR^{2n}$ under
the standard identification of  $\CC$ with $\RR^2$. 

We extend these
notions to complex analytic manifolds and their subsets in the
obvious way:\\
For a  complex analytic manifold $M$ and  $S\subseteq M$ we say that $S$ is {\em a subanalytic subset of $M$} if $S$ is a subanalytic subset of $M$ considered as a real analytic manifold.  For such $S$ we will denote by $\rdim(S)$ its dimension (as a subanalytic subset of $M$).

Recall that for a complex analytic manifold $M$, a subset $S\subseteq M$ is called {\em locally complex analytic 
subset of $M$} if for every $p\in S$ there is open $V_p\subseteq M$ containing $p$ such that $S\cap V_p$ is the zero locus of finitely many holomorphic on $V_p$ functions.\\
A locally subanalytic subset $S$ of a complex analytic manifold $M$ is called {\em complex analytic subset of $M$} if, in addition, $S$ is closed in $M$. 
\newpage

\section{Removing singularities from complex analytic sets}

\begin{thm} Let $M$ be a complex analytic manifold,
$E\sub M$ a closed subanalytic subset of $M$ and $A$ a complex
analytic subset of $M\setminus E$ which is also a subanalytic
subset of $M$.

Assume that for every open $V\sub M$, we have $\rdim(V\cap E)\leq
\rdim(V\cap A)-2$. Then $Cl(A)$, the closure of $A$ in $M$, is a
complex analytic subset of $M$.
\end{thm}

The theorem is an immediate consequence of Shiffman's Theorem in
the case when $A$ is of pure dimension. However, it is false as
stated without the subanalyticity assumption:\\
 Take $M=\CC^3$,
$E=\{(0,y,1)\in \CC^3\}$, and $A=\{(x,e^{1/x},1):x\neq 0\}\cup
\{(0,y,z):z\neq 1\}$.

Using properties of  geometric-analytic categories one can derive
from  the above theorm the following useful result (which again
fails without the subanalyticity  assumption).
\begin{thm}[Closure Theorem] Let $M$ be a complex
manifold and $E\sub M$ a complex analytic subset of $M$ (of
arbitrary dimension). If $A$ is a complex analytic subset of
$M\setminus E$ which is also a subanalytic subset of $M$ then
$\Cl(A)$ is a complex analytic subset of $M$.
\end{thm}

Another corollary  is the following:
\begin{thm}[The Union Theorem] Let $M$ be a complex manifold and $\{A_n:n\in {\mathbb N}\}$ a family of locally complex analytic subsets  of $M$. If $A=\bigcup_{n\in
\mathbb N} A_n$ is a closed subset of $M$ which is also
subanalytic in $M$ then $A$ is a complex analytic subset of $M$.
\end{thm}

Again, note that the theorem fails without the subanalyticity
assumption, even for a union of two locally analytic sets: Take
$A_1$ be the graph of the map $e^{1/z}$ in $\CC-\{0\}\times \CC$
and $A_2=\{0\}\times \CC$.

\section{Holomorphic maps and subanalytic sets}

\subsection{The image under holomorphic maps}

\begin{thm} Let $f:M\rarr N$ be a holomorphic map between
complex analytic  manifolds, $A$ a complex analytic subset of $M$.
If $f(A)$ is a closed and subanalytic subset of $N$ then $f(A)$ is
a complex analytic subset of $N$.
\end{thm}
\noindent{\bf Remarks}
\\1.  Notice that Remmert's Proper Mapping
Theorem follows from the above theorem since, in the above notation, if
$f:M\to N $ is proper then $f(A)$ is necessarily a subanalytic
closed subset of $N$. However, the theorem still applies to cases
where $f$ is not proper and yet $f(A)$ is a subanalytic subset of
$N$.
\\2. The theorem is false without the subanalyticity assumption on
$f(A)$. An example is the projection of the set $\{(0,0)\}\cup\{ (1/n,n) : n\in \NN, n\not=0\}\subseteq \CC^2$ onto the first coordinate.

\subsection{Meromorphic maps}

The following theorem also fails without subanalyticity assumption. (An easy example is the function $e^{1/z}$.) 

\begin{thm}
Let $M, N$ be complex analytic manifolds and $E\sub M$ a  closed
subanalytic subset of $M$ such that $\rdim E\leq \rdim M-2$. If
$f:M\setminus E\rarr N$ is a holomorphic map  whose graph is a
subanalytic subset of $M\times N$ then $f$ is a  meromorphic map
from $M$ to $N$, namely, the closure of the graph of $f$ is a
complex analytic subset of $M\times N$.

If $N=\CC$ then $f$ can be written locally, at every point of $M$,
as a quotient of two holomorphic functions.
\end{thm}

\begin{cor} Let $M, N$ be complex analytic manifolds, $X\sub M$ a
complex analytic subset of $M$, $Y\sub N$ a complex analytic
subset of $N$, and $f:\Reg(X)\rarr \Reg(Y)$ a bi-holomorphism.
Then $f$ extends to a bi-meromorphism from $X$ to $Y$ if and only
if the graph of $f$ is a subanalytic subset of $M\times N$.
\end{cor}

\end{document}